\documentclass[leqno]{article}
\usepackage{bcp04e}
\usepackage{amssymb}
\usepackage{xypic}

\newdefin{definition}{Definition}[section]
\newtheorem{theorem}{Theorem}[section]
\newtheorem{proposition}[theorem]{Proposition}
\newtheorem{corollary}[theorem]{Corollary}
\newtheorem{lemma}[theorem]{Lemma}
\newdefin{remark}[theorem]{Remark}
\newdefin{example}[theorem]{Example}

\def\Mmon{\kern0.3em\rule{0.04em}{0.52em}\kern-.35em\gtrdot}
\def\mmon{\kern0.17em\lower0.1ex\hbox{\rule{0.025em}{0.43em}}\kern-.105em\gtrdot}

\begin{document}
\mathclass{Primary 46L50; Secondary 60E10.}

\abbrevauthors{U. Franz}
\abbrevtitle{Multiplicative Monotone Convolutions}

\title{Multiplicative Monotone Convolutions}

\author{Uwe Franz}
\address{Institut f\"ur Mathematik und Informatik \\
Ernst-Moritz-Arndt-Universit\"at Greifswald \\
Friedrich-Ludwig-Jahn-Stra{\ss}e 15 a  \\
D-17487 Greifswald, Germay \\
Email: franz@uni-greifswald.de}
\thanks{Work supported in part by the European Community's Human Potential Programme under contract HPRN-CT-2002-00279 QP-Applications and a DAAD-KBN cooperation}

\maketitlebcp

\begin{abstract}
Recently, Bercovici has introduced multiplicative convolutions based on
Muraki's monotone independence and shown that these convolution of probability
measures correspond to the composition of some function of their Cauchy
transforms. We provide a new proof of this fact based on the combinatorics of
moments. We also give a new characterisation of the probability measures that
can be embedded into continuous monotone convolution semigroups of probability
measures on the unit circle and briefly discuss a relation to Galton-Watson
processes.
\end{abstract}

\section{Introduction.}
In quantum probability there exist several natural notions of independence,
see \cite{muraki03} and the references therein. These allow to define new
convolutions for probability measures, cf.\
\cite{voiculescu+dykema+nica92,voiculescu97,speicher+woroudi93,muraki00}.

Bercovici \cite{bercovici04} defined multiplicative monotone convolutions for
probability measures on the unit circle and on the half line. He showed that
with an appropriate function of the Cauchy transform these multiplicative
convolutions can be calculated by composition of those functions, similar to
Muraki's result \cite[Theorem 3.1]{muraki00} for the additive monotone
convolution. In this paper we give a new proof of Bercovici's result based on
the combinatorics of moments, see Theorem \ref{thm-operators}. Using Berkson
and Porta's \cite{berkson+porta78} characterization of composition semigroups, one can deduce a
characterization of continuous convolution semigroups for the monotone
convolution, see \cite[Theorem 4.6]{bercovici04} or Theorem
\ref{thm-levy-khintchine-circle} for the case of probability measures on the
unit circle.

This paper is organized as follows.

In Section \ref{sec-mon} we recall the definition of monotone independence and
the monotone product of algebraic and quantum probability spaces. In Section
\ref{sec-mon-cond} we show that the monotone product is actually a special
case of the conditionally free product introduced in
\cite{bozejko+speicher91b,bozejko+leinert+speicher96}.

Sections \ref{sec-op}, \ref{sec-conv}, and \ref{sec-levy} contain the main
results on the multiplicative monotone convolution. We formulate a slightly
modified version of a theorem by Bercovici that shows that these convolutions
can be calculated by taking the composition of appropriate functions of the
Cauchy transform of the measures, see Theorem \ref{thm-operators} and
Corollaries \ref{cor-unitary} and \ref{cor-pos}. We also state a
L\'evy-Khintchine type characterization of all continuous convolution
semigroups for the monotone convolution of probability measures on the unit
circle, see Theorem \ref{thm-levy-khintchine-circle}.

In Section \ref{sec-galton}, we show that the problem of embedding a
probability measure on the unit circle into a continuous monotone convolution
semigroup is very similar to the problem of embedding a discrete-time
Markovian branching process (or Galton-Watson process) into a continuous-time
Markovian branching process. In Section \ref{sec-embed} we adapt a
characterization of embeddable branching processes due to Gorya\u{\i}nov
\cite{goryainov93} to our situation.

Finally, in the Appendix we discuss the multiplicative monotone convolution of
probability measures on the half line and show that there exist two natural, but
inequivalent definitions. One of them is equivalent to the definition due to
Bercovici and can be treated by similar methods as the
multiplicative monotone convolution of measures on the unit circle., cf.\
\cite{bercovici04}.

\section{Monotone Independence.}\label{sec-mon}

In this section we present the definition of monotone independence and its main properties. 

By an {\em algebraic probability space} we mean a pair
$(\mathcal{A},\varphi)$ consisting of a
unital algebra $\mathcal{A}$ and a unital functional
$\varphi:\mathcal{A}\to\mathbb{C}$. Assume that we have two algebraic
probability spaces $(\mathcal{A}_1,\varphi_1)$ and
$(\mathcal{A}_2,\varphi_2)$, such that the first algebra has a
decomposition $\mathcal{A}_1=\mathbb{C}\mathbf{1}\oplus\mathcal{A}_1^0$ (direct
sum as vector spaces), where $\mathcal{A}_1^0$ is a subalgebra of
$\mathcal{A}_1$. Then we define the algebraic monotone product
$(\mathcal{A},\varphi)$ of $(\mathcal{A}_1,\varphi_1)$ and
$(\mathcal{A}_2,\varphi_2)$ as follows, see also \cite{muraki01,muraki03}. The algebra $\mathcal{A}=\mathcal{A}_1\coprod
\mathcal{A}_2$ is the free product of $\mathcal{A}_1$ and
$\mathcal{A}_2$ with identification of the units of $\mathcal{A}_1$ and
$\mathcal{A}_2$. The unital functional
$\varphi=\varphi_1\triangleright\varphi_2:\mathcal{A}\to\mathbb{C}$ is determined by the condition
\begin{equation}\label{def-alg-mon}
\varphi(b_1a_1b_2\cdots a_{n-1}b_n)=\varphi_1(a_1\cdots
a_{n-1})\varphi_2(b_1)\cdots \varphi_2(b_n)
\end{equation}
for $n\in\mathbb{N}$ and all $a_1,\ldots,a_{n-1}\in\mathcal{A}_1^0$,
$b_1,\ldots,b_n\in\mathcal{A}_2$.

Let now $\mathcal{A}_1,\mathcal{A}_2\subseteq\mathcal{B}$ be two such
algebras, which are contained in an algebraic probability space
$(\mathcal{B},\Phi)$ and denote by
$j_1:\mathcal{A}_1\to\mathcal{B}$, $j_2:\mathcal{A}_2\to\mathcal{B}$
the inclusion maps. Then the universal property of the free
product of algebras implies that there exists a unique homomorphism
$j:\mathcal{A}_1\coprod\mathcal{A}_2\to B$ such that the following
diagram commutes
\[
\xymatrix{
& \mathcal{B} & \\
\mathcal{A}_1\ar[ur]^{j_1}\ar[r]_{i_1} &
\mathcal{A}_1\coprod\mathcal{A}_2 \ar[u]|-j & \mathcal{A}_2\ar[ul]_{j_2}\ar[l]^{i_2}
}
\]
where are $i_1:\mathcal{A}_1\to\mathcal{A}_1\coprod\mathcal{A}_2$ and
$i_2:\mathcal{A}_2\to\mathcal{A}_1\coprod\mathcal{A}_2$ are the
canonical inclusion maps.

The subalgebras $\mathcal{A}_1,\mathcal{A}_2$ are called {\em monotonically
independent} w.r.t.\ $\Phi$, if
\[
\Phi\circ j = (\Phi\circ j_1)\triangleright (\Phi\circ j_2)
\]
cf.\ \cite{franz02}

We will call a triple $(\mathcal{A},\mathcal{H},\Omega)$ consisting of
a Hilbert space $\mathcal{H}$, a unit vector $\Omega\in\mathcal{H}$,
and a subalgebra $\mathcal{A}\subseteq\mathcal{B}(\mathcal{H})$ a {\em
  quantum probability space}.

If we have an algebraic probability space $(\mathcal{A},\varphi)$, whose
algebra has an involution such that $\Phi$ is even a state, and if for
all $a\in\mathcal{A}$ there exists a constant $C_a\ge 0$ such that the
inequality
\[
\Phi(x^*a^*ax)\le C_a\Phi(x^*x)
\]
holds for all $x\in\mathcal{A}$, then the GNS representation
$(H_\varphi,\pi_\varphi,\Omega_\varphi)$ of $(\mathcal{A},\Phi)$
yields a quantum probability space
$(\pi_\varphi(\mathcal{A}),H_\varphi,\Omega_\varphi)$. If two
subalgebras
$\mathcal{A}_1=\mathbb{C}\mathbf{1}\otimes\mathcal{A}_1^0,\mathcal{A}_2\subseteq\mathcal{A}$
are monotonically independent in $(\mathcal{A},\varphi)$, then
$\pi_\varphi(\mathcal{A}_1^0)$ and $\pi_\varphi(\mathcal{A}_2)$ are
monotonically independent in
$(\pi_\varphi(\mathcal{A}),H_\varphi,\Omega_\varphi)$ in the sense of
the following definition.

\begin{definition}\label{def-mon-indep}
Let $\mathcal{H}$ be a Hilbert space, $\Omega\in\mathcal{H}$ a unit vector, and define a state $\Phi:\mathcal{B}(\mathcal{H})\to\mathbb{C}$ on the algebra of bounded operators on $\mathcal{H}$ by
\[
\Phi(X)=\langle \Omega, X\Omega\rangle, \qquad\mbox{ for } X\in\mathcal{B}(\mathcal{H}).
\]
Two subalgebras $\mathcal{A}_1,\mathcal{A}_2\subseteq\mathcal{B}(\mathcal{H})$ are called {\em monotonically independent} w.r.t.\ $\Omega$, if the following two conditions are satisfied.
\begin{description}
\item[(a)]
For all $X,Z\in\mathcal{A}_1$, $Y\in\mathcal{A}_2$, we have
\[
XYZ = \Phi(Y)XZ.
\]
\item[(b)]
For all $Y\in\mathcal{A}_1$, $X,Z\in\mathcal{A}_{2}$,
\[
\Phi(XYZ)=\Phi(X)\Phi(Y)\Phi(Z).
\] 
\end{description}
Two operators $X,Y\in\mathcal{B}(\mathcal{H})$ are called
monotonically independent w.r.t.\ $\Omega$, if the subalgebras
$\mathcal{A}_1={\rm alg}(X)={\rm span}\{X^k|k=1,2,\ldots\}$ and $\mathcal{A}_2={\rm alg}(Y)={\rm span}\{Y^k|k=1,2,\ldots\}$ are monotonically independent.
\end{definition}

\begin{proposition}\label{prop-mon-prod}
Let $(\mathcal{A}_i,\mathcal{H}_i,\Omega_i)$, $i=1,2$, be two quantum probability spaces, and denote the states associated to $\Omega_1$ and $\Omega_2$ by $\Phi_1$ and $\Phi_2$, respectively.

Then there exists a quantum probability space $(\mathcal{A},\mathcal{H},\Omega)$ and two injective state-preser\-ving homomorphisms $J_i:\mathcal{A}_i\to\mathcal{A}$, $i=1,2$, such that the images $J_1(\mathcal{A}_1)$ and $J_2(\mathcal{A}_2)$ are monotonically independent w.r.t.\ $\Omega$.
\end{proposition}
\Proof
We set $\mathcal{H}=\mathcal{H}_1\otimes\mathcal{H}_2$ and $\Omega=\Omega_1\otimes\Omega_2$. Denote by $P_2$ the orthogonal projection on $\mathbb{C}\Omega_2\subseteq\mathcal{H}_2$.

We define the embeddings $J_i:\mathcal{A}_i\to\mathcal{B}(\mathcal{H})$ by
\begin{eqnarray*}
J_1(X) &=& X\otimes P_2, \qquad \mbox{ for } X\in\mathcal{A}_1, \\
J_2(X) &=& \mathbf{1}\otimes X, \qquad \mbox{ for } X\in\mathcal{A}_2.
\end{eqnarray*}
For $\mathcal{A}$ we take the subalgebra generated by $J_1(\mathcal{A}_1)$ and $J_2(\mathcal{A}_2)$. It is clear that $J_1$ and $J_2$ are injective, state-preserving homomorphisms.

A simple calculation shows that $J_1(\mathcal{A}_1)$ and $J_2(\mathcal{A}_2)$ are monotonically independent w.r.t.\ $\Omega$. E.g., for products of the form $J_1(X_1)J_2(Y)J_1(X_2)$, $X_1,X_2\in\mathcal{A}_1$, $Y\in\mathcal{A}_2$, we get
\begin{eqnarray*}
J_1(X_1)J_2(Y)J_1(X_2) &=& (X_1\otimes P_2)(\mathbf{1}\otimes Y)(X_1\otimes P_2) = (X_1X_2)\otimes P_2YP_2 \\
&=& \Phi\big(J_2(Y)\big) J_1(X_1)J_1(X_2).
\end{eqnarray*}
On the other hand, for $J_2(Y_1)J_1(X)J_2(Y_2)$, $X\in\mathcal{A}_1$, $Y_1,Y_2\in\mathcal{A}_2$, we get
\begin{eqnarray*}
\Phi\big(J_2(Y_1)J_1(X)J_2(Y_2)\big) &=& \langle\Omega_1\otimes\Omega_2,(\mathbf{1}\otimes Y_1)(X\otimes P_2)(\mathbf{1}\otimes Y_2)\Omega_1\otimes\Omega_2\rangle \\
&=& \langle\Omega_1\otimes\Omega_2,X\otimes (Y_1PY_2)\Omega_1\otimes\Omega_2\rangle \\
&=& \Phi_1(X)\Phi_2(Y_1)\Phi_2(Y_2) =\Phi\big(J_2(Y_1)\big)\Phi\big(J_1(X)\big)\Phi\big(J_2(Y_2)\big).
\end{eqnarray*}
\endproof

We will call the quantum probability space  $(\mathcal{A},\mathcal{H},\Omega)$ constructed in the previous proposition the {\em monotone product} of $(\mathcal{A}_1,\mathcal{H}_1,\Omega_1)$ and $(\mathcal{A}_2,\mathcal{H}_2,\Omega_2)$. When there is no danger of confusion, we shall identify the algebras $\mathcal{A}_1$ and $\mathcal{A}_2$ with their images $J_1(\mathcal{A}_1)$ and $J_2(\mathcal{A}_2)$, respectively.

The monotone product is associative and can be extended to more than two factors, see \cite{franz01}. But it is not commutative.

The embedding $J_1:\mathcal{A}_1\to\mathcal{A}$ is not unital and the product is not trace-preserving. If $\Phi_1|_{\mathcal{A}_1}$ is not identically equal to zero, then the calculation
\[
\Phi_1(X)\Phi_2(Y_1Y_2)=\Phi(XY_1Y_2)=\Phi(Y_2XY_1)=\Phi_1(X)\Phi_2(Y_1)\Phi_2(Y_2)
\]
for all $X\in\mathcal{A}_1$, $Y_1,Y_2\in\mathcal{A}_2$ shows that $\Phi$ can only be a trace on $\mathcal{A}$, if $\Phi_2|_{\mathcal{A}_2}$ is a homomorphism.

\section{Relation of monotone independence and conditional free
  independence.}\label{sec-mon-cond}

We recall now the definition of the conditional free product of algebraic
probability spaces and show that the monotone product is contained as a
special case.

Let $(\mathcal{A}_1,\varphi_1,\psi_1)$ and
$(\mathcal{A}_2,\varphi_2,\psi_2)$ be two unital algebras, equipped
with two unital functionals. Recall that the conditionally free
product\cite{bozejko+speicher91b,bozejko+leinert+speicher96} of $(\mathcal{A}_1,\varphi_1,\psi_1)$ and
$(\mathcal{A}_2,\varphi_2,\psi_2)$ is defined as the triple
$(\mathcal{A},\varphi,\psi)$, where $\mathcal{A}=\mathcal{A}_1\coprod
\mathcal{A}_2$ is the free product of $\mathcal{A}_1$ and
$\mathcal{A}_2$ with identification the units of $\mathcal{A}_1$ and
$\mathcal{A}_2$. The unital functionals $\varphi$ and $\psi$ on $\mathcal{A}=\mathcal{A}_1\coprod
\mathcal{A}_2$ can be defined by the conditions
\begin{equation}\label{def-cond-free}
\varphi(a_1a_2\cdots a_n)=\varphi_{\epsilon(1)}(a_1)\cdots\varphi_{\epsilon(n)}(a_n) \quad \mbox{
  and }\quad \psi(a_1a_2\cdots a_n)=0
\end{equation}
for all $n\in \mathbb{N}$ and all $a_i\in\mathcal{A}_{\epsilon(i)}$
with $\epsilon(i)\in\{1,2\}$,
$\epsilon(1)\not=\epsilon(2)\not=\cdots\not=\epsilon(n)$ and
$\psi_{\epsilon(1)}(a_1)=\cdots=\psi_{\epsilon(n)}(a_n)=0$. The
functional $\psi$ is simply the free product $\psi_1*\psi_2$ of $\psi_1$ and $\psi_2$,
cf.\ \cite{voiculescu+dykema+nica92,voiculescu97}. We will denote
$\varphi$ by
\[
\varphi=\varphi_1\,{}_{\psi_1}\kern-.5em*\kern-.3em{}_{\psi_2}\,\varphi_2.
\]
The product defined in this way for triples
$(\mathcal{A},\varphi,\psi)$ can be shown to be commutative and
associative, cf. \cite{bozejko+speicher91b,bozejko+leinert+speicher96}.

Taking pairs of the form $(\mathcal{A}_1,\varphi_1,\varphi_1)$ and
$(\mathcal{A}_2,\varphi_2,\varphi_2)$, one obtains the free
product also for the first functional, i.e.
\[
\varphi_1\,{}_{\varphi_1}\kern-.5em*\kern-.3em{}_{\varphi_2}\,\varphi_2
= \varphi_1*\varphi_2.
\]
Suppose now that the algebras $\mathcal{A}_1$ and $\mathcal{A}_2$
have decompositions
$\mathcal{A}_i=\mathbb{C}\mathbf{1}\oplus\mathcal{A}_i^0$, $i=1,2$, as
a direct sum of vector spaces, such that the $\mathcal{A}_i^0$ are
even subalgebras. If one defines functionals
$\delta_i:\mathcal{A}_i\to\mathbb{C}$ by
\begin{equation}\label{delta}
\delta_i(\lambda\mathbf{1}+a_0)=\lambda
\end{equation}
for $\lambda\in\mathbb{C}$, $a_0\in\mathcal{A}_i^0$, $i=1,2$, then one
obtains the boolean product
\[
\varphi_1\,{}_{\delta_1}\kern-.5em*\kern-.3em{}_{\delta_2}\,\varphi_2=\varphi_1\diamond\varphi_2,
\]
cf.\ \cite{speicher+woroudi93,bozejko+leinert+speicher96}.

Since the conditionally free product of triples of the form
$(\mathcal{A},\varphi,\delta)$ can be shown to be again of the same
form, the commutativity and associativity of the boolean product
follow immediately from this construction.

One can also obtain the monotone product from the
conditionally free product.
\begin{proposition}
Let $(\mathcal{A}_1,\varphi_1)$ and $(\mathcal{A}_2,\varphi_2)$ be two
algebraic quantum probability spaces and assume $\mathcal{A}_1$ has a
decomposition $\mathcal{A}_1=\mathbb{C}\mathbf{1}\oplus\mathcal{A}_1^0$, where
$\mathcal{A}_1^0$ is a subalgebra of $\mathcal{A}_1$. Define a unital
functional $\delta_1:\mathcal{A}_1\to\mathbb{C}$ as in Equation
(\ref{delta}).

Then we have
\[
\varphi_1\triangleright\varphi_2=\varphi_1\,{}_{\delta_1}\kern-.5em*\kern-.3em{}_{\varphi_2}\,\varphi_2
\]
\end{proposition}
\Proof
Let $n\in\mathbb{N}$, $\epsilon(1),\ldots,\epsilon(n)\in\{1,2\}$ such
that $\epsilon(1)\not=\epsilon(2)\not=\cdots\not=\epsilon(n)$, and
$a_1\in\mathcal{A}_{\epsilon(1)},\cdots,a_n\in\mathcal{A}_{\epsilon(n)}$
such that $\delta_1(a_k)=0$ if $\epsilon(k)=1$ and $\varphi_2(a_k)=0$
if $\epsilon(k)=2$. This implies $a_k\in\mathcal{A}_1^0$ for
$\epsilon(k)=1$ and therefore by Equation (\ref{def-alg-mon})
\[
\varphi_1\triangleright\varphi_2(a_1a_2\cdots
a_n)=\prod_{k:\epsilon(k)=2} \varphi_2(a_k)=0
\]
(If the product $a_1a_2\cdots a_n$ does not begin or end with an
element of $\mathcal{A}_2$, add $\mathbf{1}\in\mathcal{A}_2$ in order
to apply Equation (\ref{def-alg-mon})).

Therefore $\varphi_1\triangleright\varphi_2$ satisfies condition
(\ref{def-cond-free}) that defines the conditionally free product
$\varphi_1\,{}_{\delta_1}\kern-.5em*\kern-.3em{}_{\varphi_2}\,\varphi_2$.
\endproof

With this observation, Muraki's formula \cite[Theorem 3.1]{muraki00} for the
additive monotone convolution can be deduced from the analytic theory of the
additive conditionally free convolution developed in \cite{bozejko+leinert+speicher96}.

\section{Products of monotonically independent operators.}\label{sec-op}

For a bounded operator $X$ in a quantum probability space
$(\mathcal{B}(\mathcal{H}),\mathcal{H},\Omega)$ we define
\[
\psi_X(z)=\left\langle \Omega, \frac{zX}{1-zX} \Omega\right\rangle
\]
and
\[
K_X(z)= \frac{\psi_X(z)}{1+\psi_X(z)}
\]
for $|z|<1/||X||$.

The following theorem is similar to \cite[Theorem 2.2]{bercovici04}. Below we
provide a new proof.

\begin{theorem}\label{thm-operators}
Let $(B(\mathcal{H}),\mathcal{H},\Omega)$ be a quantum probability
space and $\mathcal{A}_1,\mathcal{A}_2 \subseteq B(\mathcal{H})$ two
monotonically independent subalgebras. Let $V_1,V_2\in
\mathbb{C}\mathbf{1}+\mathcal{A}_1$, such that $V_2V_1-\mathbf{1}\in
\mathcal{A}_1$ and $W\in\mathcal{A}_2$.
Then we have
\[
K_{V_1WV_2}(z)=K_{V_1V_2}\big(K_W(z)\big)
\]
for all $|z|<\min(1/||V_1WV_2||,1/||W||)$.
\end{theorem}
\Proof
Let $M=\max\big(||V_1WV_2||,||W||(||V_1V_2||+2)\big)$ and $|z|<1/M$. Then we have
\begin{eqnarray*}
\frac{zV_1WV_2}{1-zV_1WV_2} &=& \sum_{n=1}^\infty (zV_1WV_2)^n = \sum_{n=1}^\infty z^n V_1 \underbrace{W (X+\mathbf{1}) W \cdots
  W(X+\mathbf{1})} W V_2 \\
&&\hspace{55mm}n-1 \mbox{ times} \\
&=& \sum_{n=1}^\infty z^n \sum_{k=1}^n \sum_{{\nu_1,\ldots,\nu_k\ge 1}\atop{\nu_1+\cdots+\nu_k=n}} V_1
W^{\nu_1} X W^{\nu_2} X \cdots X W^{\nu_k} V_2,
\end{eqnarray*}
where $X=V_2V_1-\mathbf{1}$.

Using properties (a) and (b) in Definition \ref{def-mon-indep}, we get
\begin{eqnarray*}
\psi_{V_1WV_2}(z) &=& \left\langle\Omega,
\frac{zV_1WV_2}{1-zV_1WV_2}\Omega\right\rangle \\
&=& \sum_{n=1}^\infty z^n \sum_{k=1}^n \sum_{{\nu_1,\ldots,\nu_k\ge
    1}\atop{\nu_1+\cdots+\nu_k=n}}\left\langle\Omega,V_1X^{k-1}V_2\Omega\right\rangle\left\langle\Omega,W^{\nu_1}\Omega\right\rangle\cdots\left\langle\Omega,W^{\nu_k}\Omega\right\rangle \\
&=& \sum_{k=1}^\infty
\left\langle\Omega,V_1(V_2V_1-\mathbf{1})^{k-1}V_2\Omega\right\rangle
\big(\psi_W(z)\big)^k \\
&=& \sum_{k=1}^\infty
\left\langle\Omega,V_1V_2(V_1V_2-\mathbf{1})^{k-1}\Omega\right\rangle
\big(\psi_W(z)\big)^k \\
&=&
\sum_{k=1}^\infty\left\langle\Omega,\psi_W(z)V_1V_2\frac{1}{\mathbf{1}-\psi_W(z)(V_1V_2-\mathbf{1})}\Omega\right\rangle
\\
&=&\sum_{k=1}^\infty\left\langle\Omega,\frac{\frac{\psi_W(z)}{1+\psi_W(z)}V_1V_2}{\mathbf{1}-\frac{\psi_W(z)}{1+\psi_W(z)}V_1V_2}\Omega\right\rangle
= \psi_{V_1V_2}\big(K_W(z)\big).
\end{eqnarray*}
By uniqueness of analytic continuation, we get
\[
K_{V_1WV_2}(z)=K_{V_1V_2}\big(K_W(z)\big)
\]
for all $|z|<\min(1/||V_1WV_2||,1/||W||)$.
\endproof

\begin{corollary}\label{cor-unitary}
Let $U,V$ be two unitary operators such that $U-\mathbf{1}$ and $V$
are monotonically independent with respect to $\Omega$. Then we have
\[
K_{UV}(z)=K_{VU}(z)=K_U\big(K_V(z)\big)
\]
for all $|z|\in\mathbb{D}=\{z\in\mathbb{C}:|z|<1\}$.
\end{corollary}

\begin{corollary}\label{cor-pos}
Let $X,Y$ be two positive operators such that $X-\mathbf{1}$ and $Y$
are monotonically independent with respect to $\Omega$. Then we have
\[
K_{\sqrt{X}Y\sqrt{X}}(z)=K_X\big(K_Y(z)\big)
\]
for all $|z|<\min(1/||\sqrt{X}Y\sqrt{X}||,1/||Y||)$.
\end{corollary}

\section{Multiplicative monotone convolution for probability measures on the
  unit circle.}\label{sec-conv}

For a probability measure $\mu$ on $S^1$ we define
\[
\psi_\mu(z)=\int_{S^1} \frac{zx}{1-zx}{\rm d}\mu(x) \quad \mbox{ and }\quad K_\mu(z)=\frac{\psi_\mu(z)}{1+\psi_\mu(z)}
\]
for $z\in\mathbb{D}=\{z\in\mathbb{C}:|z|<1\}$.

We will call $K_\mu$ the {\em K-transform} of $\mu$, it characterizes the
measure $\mu$ completely. Furthermore, for a holomorphic function
$K:\mathbb{D}\to\mathbb{D}$ there exists a probability measure $\mu$ on the
unit circle $S^1$ such that $K=K_\mu$ if and only if $K(0)=0$. This follows
from the Herglotz representation theorem, the proof is similar to
\cite[Proposition 3.3]{franz04}.

It is clear that the composition of two K-transforms is again a K-transform of
some probability measure on the unit circle. In view of Corollary
\ref{cor-unitary} this suggests the following definition.

\begin{definition}
Let $\mu,\nu$ be two probability measures on $S^1$, with transforms $K_\mu$ and
$K_\nu$. Then the unique probability measure $\mu\Mmon\nu$ on $S^1$ with
\[
K_{\mu\mmon\nu}=K_\mu\circ K_\nu
\]
is called the {\em monotone convolution} of $\mu$ and $\nu$.
\end{definition}

\begin{remark}
\begin{enumerate}
\item
The monotone convolution is weakly continuous.
\item
The monotone convolution is associative,
i.e.
\[
(\lambda\Mmon\mu)\Mmon\nu=\lambda\Mmon(\mu\Mmon\nu)
\]
for all
$\lambda,\mu,\nu$, but not commutative, i.e., in general
$\mu\Mmon\nu\not=\nu\Mmon\mu$.
\item
The Dirac measure $\delta_1$ at $1$ is a two-sided unit,
$\delta_1\Mmon\mu=\mu\Mmon\delta_1=\mu$ for all $\mu$. Right convolution by a
Dirac measure $\delta_x$ acts
as translation, i.e.\ $\mu\Mmon\delta_x=T_x\mu$, where $T_x:S^1\to S^1$ is
defined by $T_x(y)=xy$ for $x\in S^1$. But $\delta_x\Mmon\mu\not=T_x\mu$ in
general.
\item
The monotone convolution is affine in the first argument. Together with weak
continuity this implies the
following formula
\[
\mu\Mmon\nu = \int_{S^1} {\rm d}\mu(x) \delta_x\Mmon\nu.
\]
\end{enumerate}
\end{remark}

\section{L\'evy-Khintchine formula for monotone convolution semigroups.}\label{sec-levy}

We call a weakly continuous one-parameter family $(\mu_t)_{t\ge0}$ of
probability measures on the unit circle a
continuous monotone convolution semigroup, if
\[
\mu_0=\delta_1 \qquad \mbox{ and } \qquad \mu_s\Mmon\mu_t=\mu_{s+t}
\]
for all $s,t\ge 0$. By definition a one-parameter family
$(\mu_t)_{t\ge}$ is a continuous monotone convolution semigroup if and only if
the $K$-transforms $K_t=K_{\mu_t}$, $t\ge 0$ form a continuous semigroup
w.r.t.\ to composition. The continuity of the $K$-transforms is uniform in $z$
on compact sets. Our main tool for characterizing continuous monotone
convolution semigroups will be Berkson and Porta's \cite{berkson+porta78}
characterisation of composition semigroups of holomorphic maps.

\begin{theorem}\label{thm-levy-khintchine-circle} \cite[Theorem 4.6]{bercovici04}
Let $(\mu_t)_{t\ge 0}$ be a weakly continuous family of probability measures
on the unit circle, with K-transforms $(K_t)_{t\ge 0}$. Then the following are
equivalent.
\begin{description}
\item[(a)]
$(\mu_t)_{t\ge 0}$ is a continuous monotone convolution semigroup.
\item[(b)]
$(K_t)_{t\ge 0}$ is a continuous semigroups w.r.t.\ to composition.
\item[(c)]
There exists a holomorphic function $u:\mathbb{D}\to\mathbb{C}$ with $\Re\,
u(z)\ge 0$ for $z\in\mathbb{D}$ such that
$(K_t)_{t\ge 0}$ is the (unique) solution of
\[
\frac{{\rm d}K_t(z)}{{\rm d}t} = -K_t(z)u\big(K_t(z)\big)
\]
for $z\in\mathbb{D}$ and $t\ge 0$, with initial condition $K_0(z)=z$.
\end{description}
\end{theorem}
\Proof
The equivalence between (a) and (b) follows from the definition
and the continuity properties of the monotone convolution.

The equivalence between (b) and (c) is an immediate consequence of
\cite[Theorem (3.3)]{berkson+porta78}, it suffices to identify the fixed point
at zero as the Denjoy-Wolff point of the $K_t$.
\endproof
\begin{remark}
\begin{enumerate}
\item
The function $u$ in (c) can be computed from the derivative of $(K_t)_{t\ge
  0}$ in $t=0$ by
\[
u(z) = -\frac{1}{z}\left.\frac{{\rm d}}{{\rm d}t}\right|_{t=0} K_t(z),
\]
we will call it the {\em generator} of $(K_t)_{t\ge 0}$.
\item
Such a function $u$ has a unique Herglotz representation
\[
u(z)=ib + \int_{S^1} \frac{w+z}{w-z}{\rm d}\rho(w),
\]
where $b$ is a real number and $\rho$ a finite measure on $S^1$.
\end{enumerate}
\end{remark}

\section{Relation to Galton-Watson processes.}\label{sec-galton}

A probability measure $\mu$ on the unit circle is called infinitely divisible
w.r.t.\ to the monotone convolution, if for all $n\in\mathbb{N}$ there exists
a probability measure $\mu_n$ on the unit circle such that
\begin{eqnarray*}
\mu&=&\underbrace{\mu_n\Mmon\cdots\Mmon\mu_n}. \\
&& \qquad n\mbox{ times}
\end{eqnarray*}
Bercovici has shown in \cite[Theorem 4.7]{bercovici04} that all infinitely
divisible probability measures can be embedded into a continuous monotone
convolution semigroup, i.e.\ if $\mu$ is infinitely divisible w.r.t.\ to the
monotone convolution, then there exists a continuous monotone convolution
semigroup $(\mu_t)_{t\ge 0}$ such that $\mu=\mu_t$ for some $t\ge0$. And from
the previous
section it is clear this implies that the
K-transform $K_\mu$ can be embedded into a continuous composition semigroup of K-transforms.

A similar problem has been studied in the theory of Galton-Watson processes.

Let $X_{n,k}$, $n,k=1,2,\ldots$ be independent, identically distributed random
variables with values in $\mathbb{N}$ with generating function
\[
\varphi(z)=\mathbb{E}(z^{X_{n,k}})=\sum_{m=0}^\infty p_m z^m,\quad \text{ for }z\in\mathbb{D}
\]
where $p_m=\mathbb{P}(X_{n,k}=m)$. Then the associated Galton-Watson process $(Y_n)_{n\ge 0}$ is defined by
$Y_0=1$, and
\[
Y_{n+1}=\sum_{k=1}^{Y_n} X_{n,k}, \qquad \mbox{ for } n\ge 1.
\]
This process describes the evolution of a population where after each step
each individual produces a random number of offspring according to the
probabilities $(p_m)_{m\ge 0}$.

Its generating functions form a discrete composition semigroup,
\[
\mathbb{E}(z^{Y_n})=\varphi^n(z),\qquad \text{ for }z\in\mathbb{D},\quad n\in\mathbb{N}.
\]
If $\mathbb{P}(X_{n,k}=0)=0$ (i.e.\ no individual dies without offspring), then
$\varphi(0)=0$ and $\varphi$ is the $K$-transform of a probability measure
$\mu$ on $S^1$. If $(Y_n)_{n\ge 0}$ can be embedded into a continuous-time
Markovian branching process (or equivalently, if $(\varphi^n)_{n\ge 0}$ can be embedded into a
continuous composition semigroup $(\varphi_t)_{t\ge0}$ of generating
functions), then $\mu$ is infinitely divisible for the monotone convolution
and can be embedded into a continuous monotone convolution semigroup. The problem of embedding Galton-Watson processes has been studied by Gorya\u{\i}nov\cite{goryainov93,goryainov00}.

\begin{example}
Continuous-time Markovian branching processes with extinction probability $0$ can be obtained by choosing infinitesimal offspring probabilities $\lambda_j\ge 0$ for $j\ge 2$ such
  that $\alpha=\sum_{j=2}^\infty \lambda_j<\infty$, setting
\[
v(z)=\sum_{j=2}^\infty \lambda_jz^j-\alpha z, \qquad \mbox{ for } |z|\le 1,
\]
and solving the differential equation
\[
\frac{{\rm d}}{{\rm d}t}\varphi_t(z) = v\big(\varphi_t(z)\big)
\]
with initial condition $\varphi_0(z)=z$, cf.\ \cite[Theorem 4]{goryainov93}.

A simple example is the Yule process, where $v(z)=\alpha(z^k-z)$ and
\[
\varphi_t(z) =\frac{ze^{-\alpha
    t}}{\sqrt[k-1]{1-\left(1-e^{-\alpha(k-1)t}\right)z^{k-1}\,}}, \qquad t\ge 0,
\]
for some $k\in\mathbb{N}$, $k\ge 2$. This process describes a population were
the individuals are replaced by $k$ new individuals after an exponentially
distributed random time.
\end{example}

\section{On the embedding of probability measures into continuous monotone
  convolution semigroups.}\label{sec-embed}

\cite[Theorem 6]{goryainov93} and \cite[Theorem 7]{goryainov93} characterize
probability generating functions that can be embedded into composition
semigroups of probability generating functions. In this section we give a similar characterization for K-transforms
of probability measures on the unit circle that can be embedded into continuous
monotone convolution semigroups.

Let $(K_t)_{t\ge 0}$ be a continuous composition semigroups of
K-transforms. By \cite{berkson+porta78}, $K_t$ is differentiable w.r.t.\ $t$
and satisfies the differential equation
\begin{equation}\label{eq-4}
\frac{{\rm d}}{{\rm d}t} K_t(z) = v\big(K_t(z)\big) = v(z)K'_t(z)
\end{equation}
for $t\ge 0$, $z\in\mathbb{D}$, with $v$ given by
\[
v(z)=\left.\frac{{\rm d}}{{\rm d}t}\right|_{t=0} K_t(z).
\]
This equation follows from the semigroup property $K_{s+t}=K_s\circ
K_t=K_t\circ K_s$ by differentiation w.r.t.\ $s$ at $s=0$.

By Theorem \ref{thm-levy-khintchine-circle}, the function $v$ is of the form
$v(z)=-zu(z)$, with a holomorphic function $u:\mathbb{D}\to\mathbb{C}$ such
that $\Re\,u(z)\ge 0$ for $z\in\mathbb{D}$.

We will need the following lemma.

\begin{lemma}\label{lemma-unique}
Let $u:\mathbb{D}\to\mathbb{C}$, $u\not\equiv0$, be a holomorphic function such
that $\Re\,u(z)\ge 0$ for $z\in\mathbb{D}$ and set $\beta=u(0)$, $v(z)=-zu(z)$
for $z\in\mathbb{D}$.

Then, for all $t\ge 0$, the equation
\begin{equation}\label{eq-5}
v\big(f(z)\big) = v(z)f'(z), \qquad z\in\mathbb{D},
\end{equation}
has a unique solution $f$ with $f'(0)=e^{-t\beta}$.
\end{lemma}
\Proof
The proof of this lemma is borrowed from \cite[Lemma 2]{goryainov93}.

Let $(K_t)_{t\ge 0}$ be a composition semigroup of K-transforms with generator
$u$. Then all $K_t$, $t\ge 0$ satisfy Equation (\ref{eq-5}). Furthermore, the
differential equation that the $K_t$ satisfy, implies
\[
\frac{{\rm d}}{{\rm d}t}K'_t(0)=-u(0)K'_z(0)
\]
and therefore $K'_t(0)=e^{-t\beta}$, since $K_0(z)=z$ and $K'_0(0)=1$. This proves existence.

Let now $f$ be an arbitrary solution of Equation (\ref{eq-5}) with
$f'(0)=e^{-t\beta}$. Since $v$ has no zeros inside $\mathbb{D}$ other than
$z=0$, we get $f(0)=0$ by substituting $z=0$ into Equation (\ref{eq-5}).
Differentiation Equation (\ref{eq-5}) $k$ times, we can calculate the higher
derivatives of $f$ at zero from $f'(0)=e^{-t\beta}$ and the derivatives of $v$
at zero. This proves uniqueness.
\endproof

\begin{remark}
Let $(K_t)_{t\ge0}$ be the K-transforms of a continuous monotone convolution semigroup
$(\mu_t)_{t\ge 0}$ with generator $u$. Then $K'_t(0)=e^{-tu(0)}$ is the first
moment of $\mu_t$, i.e.\
\[
e^{-tu(0)}=\int_{S^1} x {\rm d}\mu_t, \qquad \mbox{ for }t\ge 0.
\]
\end{remark}

We come to the main result of this section.
\begin{theorem}
Let $\mu$ be a probability measure on the unit circle $S^1$ that is not
concentrated in one point.

Then $\mu$ can be embedded into a continuous monotone convolution semigroup if and only
if $K'_\mu(z)\not=0$ for all $z\in\mathbb{D}$ and there exists a locally
uniform limit
\[
\lim_{n\to\infty} -\frac{K^n_\mu(z)}{(K^n_\mu)'(z)}=v(z),
\]
in $\mathbb{D}$ that is of the form $v(z)=\alpha z u(z)$ with a non-zero constant
$\alpha\in\mathbb{C}$ and a holomorphic function
$u:\mathbb{D}\to\mathbb{C}$ such that $\Re\,u(z)\ge 0$ for $z\in\mathbb{D}$
and $K'_\mu(0)=e^{-t_0u(0)}$ for some some $t_0\ge 0$.
\end{theorem}
\Proof
The proof of this theorem is similar to that of \cite[Theorem 6]{goryainov93}.

Suppose that $\mu$ can be embedded into a continuous monotone convolution semigroup. Then $K_{\mu}$ can be embedded into a
composition semigroup of K-transforms $(K_t)_{t\ge 0}$. Therefore all $K_t$ are
injective and $K'_t(z)\not=0$ for all $z\in\mathbb{D}$, $t\ge0$, cf.\
\cite{berkson+porta78}. Denote by $u$ the generator of $(K_t)_{t\ge 0}$ and
define $v$ by $v(z)=-zu(z)$ for $z\in\mathbb{D}$. By the Denjoy-Wolff theorem we get $\lim_{t\to\infty}K_t(z)=0$ and $\lim_{t\to\infty}K'_t(z)=0$
locally uniformly for all $z\in\mathbb{D}$. Therefore
\[
\lim_{t\to \infty} \frac{v\big(K_t(u)\big)}{K_t(z)}=v'(0)=-u(0).
\]
With the right-hand-side of Equation (\ref{eq-4}) this implies
\[
\lim_{n\to\infty} -\frac{K^n_\mu(z)}{(K^n_\mu)'(z)}=\lim_{t\to\infty}
-\frac{K_t(z)}{K'_t(z)}
=\lim_{t\to\infty}-\frac{K_t(z)v(z)}{v\big(K_t(z)\big)}= -\frac{v(z)}{v'(0)}=-z\frac{u(z)}{u(0)}.
\]
The limit is of the form required in the theorem with the constant $\alpha=-1/u(0)$.

To show the converse, let now $K_\mu$ be a K-transform satisfying the
conditions of the theorem with $v(z)=\alpha z u(z)$, $\alpha$ and $u$ as described in the theorem.

Let $(K_t)_{t\ge 0}$ be the composition semigroup of K-transforms with
generator $u$. Then the $K_t$ satisfy
\[
v\big(K_t(z)\big) = v(z)K'_t(z), \qquad \mbox{ for } t\ge 0, \quad z\in\mathbb{D}.
\]
The conditions of the theorem imply that $K_\mu$ is also a solution of the same
equation,
\[
v(z)=\lim_{t\to\infty} -\frac{K_\mu^{n+1}(z)}{(K_\mu^{n+1})'(z)}
=\lim_{t\to\infty} -\frac{K_\mu^n\big(K_\mu(z)\big)}{K'_\mu(z)(K_\mu^n)'(z)} =
\frac{v\big(K_\mu(z)\big)}{K'_\mu(z)}.
\]
The uniqueness in Lemma \ref{lemma-unique} now implies $K_\mu=K_{t_0}$.
\endproof
\begin{remark}
Let $\mu=\delta_x$ be concentrated in one point $x=e^{i\varphi}\in S^1$. Then
we get $\psi_{\delta_x}(z)=\frac{xz}{1-xz}$ and $K_\mu(z)=e^{i\varphi}z$ and
$\mu$ can be embedded into the continuous convolution semigroups
$(\mu_t^{(k)})_{t\ge 0}$ given by $\mu^{(k)}_t=\delta_{e^{it(\varphi+2\pi
    k)}}$, $k\in\mathbb{Z}$.
\end{remark}

\section{Appendix: Multiplicative monotone convolution for
 probability measures on $\mathbb{R}_+$.}

Just as there are many different ways to define multiplicatively a positive
operator from two given positive operators, there are different definitions of multiplicative monotone
convolutions of two probability measures $\mu$ and $\nu$ on
$\mathbb{R}_+$. Two possible choices are to take positive self-adjoint operators $X$ and $Y$, whose distributions
are given by $\mu$ and $\nu$, resp., such that $X-\mathbf{1}$ and
$Y-\mathbf{1}$ are monotonically independent, and to define the
convolution of $\mu$ and $\nu$
as the distributions of $\sqrt{X}Y\sqrt{X}$ or
$\sqrt{Y}X\sqrt{Y}$.

By Corollary \ref{cor-pos} the K-transform of $\sqrt{X}Y\sqrt{X}$ is equal to
the composition of the K-transforms of $X$ and $Y$. Therefore this gives a
definition which is equivalent to the one chosen by Bercovici, cf.\
\cite{bercovici04}. 

We will show below that choosing the distribution of $\sqrt{Y}X\sqrt{Y}$
as the convolution of the distributions of $X$ and $Y$ leads to an
inequivalent definition.

The operators  $\sqrt{X}Y\sqrt{X}$ and
$\sqrt{Y}X\sqrt{Y}$ have the same spectrum, except for $0$. More
precisely,
$\sigma(\sqrt{X}Y\sqrt{X})\backslash\{0\}=\sigma(\sqrt{Y}X\sqrt{Y})\backslash\{0\}$,
since $\sqrt{X}Y\sqrt{X}=AB$ and $\sqrt{Y}X\sqrt{Y}=BA$ with
$A=\sqrt{X}\sqrt{Y}$ and $B=\sqrt{Y}\sqrt{X}$.

But the following example shows that, unlike in the free case where one
works with tracial states, here the distributions of
$\sqrt{X}Y\sqrt{X}$ and $\sqrt{Y}X\sqrt{Y}$ are in general different
and therefore we have two different multiplicative monotone
convolutions for probability measures on $\mathbb{R}_+$.

\begin{example}
Consider the positive definite $2\times 2$-matrix
\[
M(a)=\left(\begin{array}{cc} 1 & a \\ a &
  1\end{array}\right)=\frac{1}{\sqrt{2}}\left(\begin{array}{cc} 1 & 1
  \\ 1 & -1 \end{array}\right)\left(\begin{array}{cc} 1+a & 0
  \\ 0 & 1-a \end{array}\right)\frac{1}{\sqrt{2}}\left(\begin{array}{cc} 1 & 1
  \\ 1 & -1 \end{array}\right),
\]
with $a\in(0,1)$. Then we have
\[
\left\langle\left(\begin{array}{c} 0 \\ 1 \end{array}\right), A^k
  \left(\begin{array}{c} 0 \\ 1 \end{array}\right)\right\rangle =
  \frac{1}{2}\left((1-a)^k+(1+a)^k\right)
\]
for $k\in\mathbb{N}$, i.e.\ the distribution of $A$ in the vector
state given by $\omega=\left(\begin{array}{c} 0 \\ 1
  \end{array}\right)$ is equal to
$\frac{1}{2}(\delta_{1-a}+\delta_{1+a})$.

A simple calculation yields
\begin{equation}\label{sqrt}
\sqrt{M(a)}=\frac{1}{2}\left(\begin{array}{ccc} \sqrt{1+a}+\sqrt{1-a} & & \sqrt{1+a}-\sqrt{1-a}
  \\ \sqrt{1+a}-\sqrt{1-a} & & \sqrt{1+a}+\sqrt{1-a} \end{array}\right).
\end{equation}

Let $a,b\in(0,1)$ and consider the pair of positive definite matrices
\begin{eqnarray*}
X&=&\mathbf{1}\otimes\mathbf{1}+\big(M(a)-\mathbf{1}\big)\otimes
P_\omega = \left(\begin{array}{cccc}
1 & 0 & 0 & 0 \\
0 & 1 & 0 & 0 \\
0 & 0 & 1 & a \\
0 & 0 & a & 1
\end{array}\right), \\
Y&=&\mathbf{1}\otimes M(b)
= \left(\begin{array}{cccc}
1 & 0 & b & 0 \\
0 & 1 & 0 & b \\
b & 0 & 1 & 0 \\
0 & b & 0 & 1
\end{array}\right),
\end{eqnarray*}
in $\mathcal{M}_2(\mathbb{C})\otimes
\mathcal{M}_2(\mathbb{C})\cong\mathcal{M}_4(\mathbb{C})$ where
$P_\omega$ denotes the orthogonal projection onto
$\omega=\left(\begin{array}{c} 0 \\ 1 \end{array}\right)$. With
respect to the vector state given by $\omega\otimes \omega$, $X-\mathbf{1}\otimes\mathbf{1}$ and $Y-\mathbf{1}\otimes\mathbf{1}$ are monotonically independent, with
distributions given by $\frac{1}{2}(\delta_{1-a}+\delta_{1+a})$ and
$\frac{1}{2}(\delta_{1-b}+\delta_{1+b})$, respectively.

As in Equation (\ref{sqrt}), we compute
\begin{eqnarray*}
\sqrt{X} &=& \left(\begin{array}{ccccccc}
1 & & 0 & & 0 & & 0 \\
0 & & 1 & & 0 & & 0 \\
0 & & 0 & & \frac{\sqrt{1+a}+\sqrt{1-a}}{2} & & \frac{\sqrt{1+a}-\sqrt{1-a}}{2} \\
0 & & 0 & & \frac{\sqrt{1+a}-\sqrt{1-a}}{2} & &
\frac{\sqrt{1+a}+\sqrt{1-a}}{2}\end{array}\right), \\
\sqrt{Y} &=& \frac{1}{2}\left(\begin{array}{ccccccc}
\sqrt{1+b}+\sqrt{1-b} & & 0 & & \sqrt{1+b}-\sqrt{1-b} & & 0 \\
0 & & \sqrt{1+b}+\sqrt{1-b} & & 0 & & \sqrt{1+b}-\sqrt{1-b} \\
\sqrt{1+b}-\sqrt{1-b} & & 0 & & \sqrt{1+b}+\sqrt{1-b} & & 0 \\
0 & & \sqrt{1+b}-\sqrt{1-b} & & 0 & &
\sqrt{1+b}+\sqrt{1-b}\end{array}\right).
\end{eqnarray*}
The eigenvalues of both $\sqrt{X}Y\sqrt{X}$ and $\sqrt{X}Y\sqrt{X}$ are
\begin{eqnarray*}
\lambda_1 &=& 1+\frac{a}{2}+\frac{1}{2}\sqrt{a^2+4(1+a)b^2}, \\
\lambda_2 &=& 1+\frac{a}{2}-\frac{1}{2}\sqrt{a^2+4(1+a)b^2}, \\
\lambda_3 &=& 1-\frac{a}{2}+\frac{1}{2}\sqrt{a^2+4(1-a)b^2}, \\
\lambda_4 &=& 1-\frac{a}{2}-\frac{1}{2}\sqrt{a^2+4(1-a)b^2},
\end{eqnarray*}
and therefore their distributions have the same support.  But
their distributions in the
vector state $\omega=\left(\begin{array}{c} 0 \\ 0 \\ 0 \\ 1
\end{array}\right)$ are different.
For example their second moments differ,
\begin{eqnarray*}
\langle\omega,\left(\sqrt{X}Y\sqrt{X}\right)^2\omega\rangle &=&
1+b^2+a^2, \\
\langle\omega,\left(\sqrt{Y}X\sqrt{Y}\right)^2\omega\rangle &=&
1+b^2+\frac{a^2}{2}\left(1+\sqrt{1-b^2}\right),
\end{eqnarray*}
(recall that we assumed $a\not=0$, $b\not=0$).
\end{example}

\section*{Acknowledgements.}

I presented the results of this paper, in particular Theorems \ref{thm-operators} and \ref{thm-levy-khintchine-circle} at
the conference ``Quantum Probability and Infinite Dimensional Analysis'' in
B\c{e}dlewo, Poland, in June 2004. I wish to thank Marek Bo{\.z}ejko and Janusz
Wysoczanski who indicated Reference \cite{bercovici04} to me.

I am also indebted to W.\ Hazod for bringing Gorya\u{\i}nov's work
\cite{goryainov93,goryainov00} to my attention.

\end{document}